\documentclass[11pt]{amsart}

\swapnumbers

\usepackage{amssymb,amsfonts}
\usepackage{euscript}

\usepackage[ps,matrix,arrow,curve,cmtip]{xy}

\title{Power operations for Morava $E$-theory of height $2$ at the
  prime $2$}
\author{Charles Rezk}
\date{ \today}
\address{Department of Mathematics \\
University of Illinois at Urbana-Champaign \\ 
Urbana, IL}
\email{rezk@math.uiuc.edu}

\thanks{The author was
partially supported by National Science Foundation grant
DMS-0505056.}

\numberwithin{equation}{section}

\makeatletter
  \let\c@subsection\c@equation
\makeatother

\theoremstyle{plain}   

\theoremstyle{remark}

\newtheorem{exer}[subsection]{Exercise}

\theoremstyle{plain}

\newcommand{\op}{{\operatorname{op}}}

\newcommand{\ra}{\rightarrow}

\newcommand{\xra}{\xrightarrow}

\newcommand{\xla}{\xleftarrow}

\DeclareMathOperator{\Tor}{Tor}

\newcommand{\powser}[1]{[\![#1]\!]}

\newcommand{\F}{\mathbb{F}}
\newcommand{\Z}{\mathbb{Z}}

\newcommand{\sm}{\wedge} 

\begin{document}

\begin{abstract}
Explicit calculations of the algebraic theory of power operations for
a specific Morava $E$-theory spectrum are given, without detailed proofs.
\end{abstract}

\maketitle


\newcommand{\Frob}{\operatorname{Frob}}

\newcommand{\rank}{\operatorname{rank}}

\newcommand{\hS}{\widehat{S}}

\section{Introduction}

In this note, we give an explicit description of the algebraic theory
of power operations for a particular Morava $E$-theory spectrum
associated to the deformations of a supersingular elliptic curve at
the prime $2$.
To understand the context behind the calculations, this note should be
read together with 
\cite{rezk-morava-e-theory-congruences}.  The purpose of this note is
to lay out some calcualtions (\S2), without 
giving complete proofs, but merely indicating how they arise from
algebraic geometry and topology (\S3 and \S4).  

Warning: there is difference of convention with
\cite{rezk-morava-e-theory-congruences}.  In that paper, it was
convenient to arrange things so that the interesting objects were
\emph{right $\Gamma$-modules}.  Here, because I'm more comfortable
carrying out computations this way, we arrange things so that the interesting
objects are \emph{left $\Gamma$-modules}.  To translate between the two
conventions, replace $\Gamma$ with $\Gamma^\op$.

I'd like to thank Matt Ando for many helpful conversations about power
operations, and Nick Kuhn for leading me to  the right way to find the
adem relations.  The calculations in this paper were aided at
numerous points by the \emph{Macualay 2} commutative algebra software package
\cite{M2}.

\section{Algebraic constructions}

In this section, we will work with the ground ring $R=\Z[a]$, a
polynomial ring on one generator.   

\subsection{The ring $\Gamma$}

We define an associative ring $\Gamma$ equipped with a ring
homomorphism $\eta\colon R\ra \Gamma$ as follows.  The ring $\Gamma$
is generated over $R$ by elements $Q_0$, $Q_1$, and $Q_2$, subject to
(i) 
\emph{commutation relations}, and (ii) \emph{adem relations}.  The
commutation relations state that the $Q_i$'s commute with elements of
$\Z\subset R$, and that
\begin{align*}
Q_0\,a &= a^2\,Q_0 - 2a \,Q_1  +6\,Q_2,
\\
Q_1\,a &= 3\,Q_0+a\,Q_2,
\\
Q_2\,a &= -a\,Q_0 +3\,Q_1.
\end{align*}
The adem relations are
\begin{align*}
Q_1Q_0 &= 2\,Q_2Q_1 - 2\,Q_0Q_2,
\\
Q_2Q_0 &= Q_0Q_1 + a\,Q_0Q_2 - 2\,Q_1Q_2.
\end{align*}

It follows using these relations that $\Gamma$ has an \emph{admissible
  basis}; that is, 
$\Gamma$ is free as a left $R$-module on the elements of the form
\[
Q_0^{j} Q_{k_1}\cdots Q_{k_r},\qquad j\geq0,\, k_i\in\{1,2\},\, r\geq0.
\]
Note that if we write $\Gamma[k]$ for the degree $k$ part of $\Gamma$,
then $\rank \Gamma[k]= 1+2+\cdots +2^k$.

\subsection{$\Gamma$-modules}

By \emph{$\Gamma$-module}, we will mean a \emph{left} $\Gamma$-module, unless
otherwise specified.  Any $\Gamma$-module is automatically an
$R$-module, via the ring homomorphism $\eta\colon R\ra \Gamma$.  Given
an $R$-module $M$, we say that a $\Gamma$-module structure on $M$
\emph{extends} the given $R$-module structure if these module
structures are compatible with respect to $\eta$.

When there is possibility of confusion, for $M$ a $\Gamma$-module,
$\gamma\in \Gamma$, $m\in M$, we write $\gamma\cdot m\in M$ for the
action of $\gamma$ on $m$.  Otherwise, we just write $\gamma m$ for
this element.

There is a standard $\Gamma$-module structure on $R$, which is
characterized as the unique one extending the $R$-module strucure on
$R$, and such that
\begin{align*}
Q_0 \cdot 1 &= 1,
\\
Q_1 \cdot 1 &= 0,
\\
Q_2 \cdot 1 &= 0.
\end{align*}

Given $\Gamma$-modules $M$ and $N$, we define their tensor product
$M\otimes N$ as follows.  The underlying $R$-module is $M\otimes_R
N$.  The $\Gamma$-module structure is given by
\begin{align*}
Q_0(x\otimes y) &= Q_0x\otimes Q_0y + 2\,Q_1x\otimes
Q_2y+2\,Q_2x\otimes Q_1y,
\\
Q_1(x\otimes y) & = Q_0x\otimes Q_1y + Q_1x\otimes Q_0y + a\,
Q_1x\otimes Q_2y +a\, Q_2x\otimes Q_1y +2\, Q_2x\otimes Q_2y,
\\
Q_2(x\otimes y) &= Q_0x\otimes Q_2y+Q_2x\otimes Q_0y+Q_1x\otimes Q_1y
+ a\, Q_2x\otimes Q_2y.
\end{align*}

The category of $\Gamma$-modules, equipped with the tensor product
$\otimes$ and unit object $R$, is a symmetric monoidal category.  

\subsection{The element $\Psi$}

Let $\Psi\in \Gamma$ be the element defined by 
\[
\Psi = Q_0Q_0 +a\,Q_0Q_1 -2\,Q_1Q_1  +a^2\,Q_0Q_2  -2a\,Q_1Q_2
+4\,Q_2Q_2.
\]
It is ``easy'' to check that $\Psi$ is in the center of
$\Gamma$, and that for a tensor product $M\otimes N$ of
$\Gamma$-modules we have $\Psi(x\otimes y)=\Psi x\otimes \Psi y$.  

\subsection{The module $\omega$}

The $\Gamma$-module $\omega$ is defined as follows.  The
underlying $R$-module of $\omega$ is a free $R$-module on one
generator $u$.  The $\Gamma$-module structure is the unique one
extending this $R$-module structure, and with
\begin{align*}
Q_0\cdot u &= 0, \\
Q_1\cdot u &= -u, \\
Q_2\cdot u &= 0.
\end{align*}
Observe that $\Psi \cdot u = -2u$.

We write $\omega^n$ for the $n$th tensor power of $\omega$, $n\geq0$.

\subsection{$\Gamma$-rings}

A \emph{$\Gamma$-ring} is a commutative monoid object in $\Gamma$-modules.
Equivalently, a $\Gamma$-ring is a commutative $R$-algebra equipped
with a $\Gamma$-module structure compatible with the $R$-module
structure, and which satisfies the \emph{cartan formulas}
\begin{align*}
Q_0(x y) &= Q_0x\,Q_0y + 2\,Q_1x\,
Q_2y+2\,Q_2x\, Q_1y,
\\
Q_1(x y) & = Q_0x\, Q_1y + Q_1x\, Q_0y + a\,
Q_1x\, Q_2y +a\, Q_2x\, Q_1y +2\, Q_2x\, Q_2y,
\\
Q_2(x y) &= Q_0x\, Q_2y+Q_2x\, Q_0y+Q_1x\, Q_1y
+ a\, Q_2x\, Q_2y.
\end{align*}

The ring $R$, with its standard $\Gamma$-module structure, is a
$\Gamma$-ring.

\subsection{Frobenius congruence and amplified $\Gamma$-rings}

We say that a $\Gamma$-ring $A$ satisfies the \emph{frobenius
  congruence} if for all $x\in A$, we have
\[
Q_0 x \equiv x^2 \mod 2A.
\]

An \emph{amplified $\Gamma$-ring} is a $\Gamma$-ring $A$, together
with a function $\theta\colon A\ra A$ which is a ``witness'' for the
frobenius congruence.  That is, the identity
\[
Q_0 x = x^2 + 2\,\theta x
\]
holds for all $x\in A$, and furthermore a number of identities
involving $\theta$ hold, all of which are identities necessarily
satisfied in a torsion free $\Gamma$-ring by a function $\theta$
satisfying the above identity.
Among these additional identities, we have
\begin{align*}
  \theta(x+y) & = \theta x + \theta y - xy,
\\
 \theta(ax) &= a^2\,\theta x -a\,Q_1x+3\,Q_2x,
\\
 \theta(xy) &= x^2\,\theta y+y^2\, \theta x +2\,\theta x\,\theta y+
 Q_1x\, Q_2y+Q_2x\,Q_1y,
\\
  Q_1\theta(x) &= Q_2Q_1x-Q_0Q_2x -Q_0x\,Q_1x -a\,Q_1x\,Q_2x -
  (Q_2x)^2,
\\
  Q_2\theta(x) &= \theta Q_1 x + a\, \theta Q_2 x - Q_1Q_2 x - Q_0x\,Q_2x,
\end{align*}
for all $x,y\in A$.  We also have $\Psi \theta = \theta \Psi$.

Amplified $\Gamma$-rings are an example of an \emph{amplified
  plethory}, as in \cite{borger-wieland-plethystic-algebra}.

The free amplified $\Gamma$-ring on one generator $x$ is a
polynomial ring 
of the form
\[
R[ \theta^j Q_{k_1}\cdots Q_{k_r} x,\, j,r\geq0,\, k_i\in \{1,2\}].
\]

\subsection{Extension of scalars}

Consider the $R$-algebras $S=\Z[a,(a^3-27)^{-1}]$ and
$\hS=\Z_2\powser{a}$.
These both admit in a unique way the structure of amplified
$\Gamma$-ring.    Note that this is not obviously the case, since $a$
is not central in $\Gamma$.  For
instance, to show that the action of $\Gamma$ can be defined on $\hS$,
one shows that $\Gamma\cdot 2R \subseteq 2R$ and $\Gamma\cdot
a^3R\subseteq 2R+aR$, and thus the 
elements of $\Gamma$ act continuously on $R$ with respect to the adic
topology.

The quotient ring $R/(2)$ is a $\Gamma$-ring, but is not an amplified
$\Gamma$-ring.  The quotient ring $R/(2,a)$ is not even a
$\Gamma$-ring.

\subsection{Trace, norm, and logarithm operators}

Define the element $T$ in $\Gamma$ by
\[
T=3\,Q_0+2a\,Q_2.
\]
For a $\Gamma$-ring $A$, we define a function $N\colon A\ra A$ by
\begin{align*}
Nx &= (Q_0x)^3 +2a\, (Q_0x)^2Q_2x -a\, Q_0x(Q_1x)^2 + a^2\,
Q_0x(Q_2x)^2 
\\ &\qquad\qquad\qquad
- 6\, Q_0\, Q_1x \, Q_2 x 
+2\,(Q_1x)^3 
 -2a\, Q_1x(Q_2x)^2
+4\,(Q_2x)^3.
\end{align*}
The function $N$ is multiplicative: $N(xy)=(Nx)(Ny)$ for any $x,y\in
A$.  

For $m\in \Z\subset R$, we have $N(m)=m^3$.  One can compute
that
\[
N(a-3) = -(a-3)^3,\qquad N(a^3-27) = -(a^3-27)^3.
\]

The operator $T$ is the ``linearization'' of $N$.  That is, if $A\ra
A/I$ is a homomorphism of $\Gamma$-rings, then for
$x\in I$ we have
\[ 
N(1+x) \equiv Tx \mod I^2.
\]

It is straightforward to show that if $A$ is an amplified
$\Gamma$-ring, and $x\in A$, then the congruence
\[
Nx \equiv x^2\, \Psi x \mod 2A
\]
holds.
Let $A^\times\subset A$ denote the group of units in $A$.  The above
congruence implies that there is a function $M\colon A^\times \ra A$
defined by
\[
1+2\,Mx = \frac{x^2\,\Psi x}{Nx}.
\]
Now define a homomorphism $\ell \colon A^\times \ra A^{\sm}_2$ by
\[
\ell(x) = \frac{1}{2} \log\left( \frac{x^2\,\Psi x}{Nx} \right) =
\sum_{k\geq1} (-1)^{k-1}\frac{2^{k-1}}{k} (Mx)^k.
\]
As an example, observe that for $D=a^3-27$ in $S=\Z[a][D^{-1}]$, we have 
\[
\ell(D) = \frac{1}{2}\log \left( \frac{D^2\, \Psi D}{N D} \right) =
\frac{1}{2} \log ( -1 ) = 0.
\]
In fact, $\ell(x)=0$ for all $x\in S^\times$.  (This is not true for
all elements of $\hS^\times$.)

\begin{exer}
Calculate the kernel of $\ell\colon \hS^\times \ra \hS$.
\end{exer}

\subsection{A Koszul complex for $\Gamma$}

Let $C_1=Rq_0+Rq_1+Rq_2$, a free left $R$-module on $3$ generators
$q_0,q_1,q_2$.  We impose a right $R$-module structure on $C_1$
(distinct from the left $R$-module structure), by
\begin{align*}
  q_0\,a &= a^2\,q_0-2a\,q_1+6\,q_2, \\
  q_1\,a &= 3\,q_0+a\,q_2,\\
  q_2\,a &= -a\,q_0+3\,q_1.
\end{align*}
(Thus $C_1$ is evidently isomorphic to the degree $1$ part
$\Gamma[1]\subset \Gamma$, as a sub-$R$-bimodule of
$\Gamma$.)

Let $C_2=Rr_1+Rr_2$, a free left $R$-module on $2$ generators
$r_1,r_2$.  We define a right $R$-module structure on $C_2$ to
coincide with the left $R$-module structure, so that $r_1a=ar_1$ and
$r_2a=ar_2$.   (We can identify $C_2$ with the kernel of
multiplication $\Gamma[1]\otimes_R \Gamma[1]\ra \Gamma[2]$.)

Given a (left) $\Gamma$-module $M$, we define a chain complex
\[
0\ra \Gamma\otimes_R C_2\otimes_R M \xra{d_2} \Gamma\otimes_R C_1\otimes_R M
\xra{d_1} \Gamma\otimes_R M\xra{d_0} M\ra 0
\]
as follows.  The tensor products are as $R$-bimodules.  The
differentials are defined by 
\begin{align*}
  d_0(\gamma\otimes m) & = \gamma m, 
\\
d_1(\gamma\otimes q_0\otimes m) &= \gamma\otimes Q_0m-\gamma
Q_0\otimes m,
\\
d_1(\gamma\otimes q_1\otimes m) &= \gamma\otimes Q_1m -\gamma
Q_1\otimes m,
\\
d_1(\gamma\otimes q_2\otimes m) &= \gamma\otimes Q_2m -\gamma
Q_2\otimes m,
\\
d_2(\gamma\otimes r_1\otimes m) &= \gamma Q_1\otimes q_0\otimes m +
\gamma(-2Q_2)\otimes q_1\otimes m + \gamma(2Q_0)\otimes q_2\otimes m
\\
&  +\gamma\otimes q_1\otimes Q_0m + \gamma\otimes (-2q_2)\otimes
Q_1 m + \gamma\otimes (2q_0)\otimes Q_2 m,
\\
d_2(\gamma\otimes r_2\otimes m) 
&= \gamma Q_2\otimes q_0\otimes m + \gamma(-Q_0)\otimes q_1\otimes m +
\gamma(-aQ_0)\otimes q_2\otimes m + \gamma(2Q_1)\otimes q_2\otimes m
\\
& +\gamma\otimes q_2\otimes Q_0 m +\gamma\otimes (-q_0)\otimes
Q_1m +\gamma\otimes (-aq_0)\otimes Q_2m + \gamma\otimes (2q_1)\otimes Q_2m.
\end{align*}
The ideas of \cite{priddy-koszul-resolutions} apply to show that this 
complex of left $\Gamma$-modules is acyclic when $M$ is flat as an
$R$-module.

\begin{exer}
Let $I\subset \Gamma$ denote the two-sided ideal generated by
$Q_0,Q_1,Q_2$.   Regard the quotient $\Gamma/I$ as a
\emph{right} $\Gamma$-module.  
Use the Koszul complex to calculate $\Tor^{\Gamma}_q(\Gamma/I, \omega^k)$
for various values of $k\geq0$.  
\end{exer}

\section{An elliptic curve}

Let $R=\Z[a]$.  Consider the curve $C$ over $R$
defined by the projective equation 
\[
Y^2Z+aXYZ+YZ^2-X^3=0
\]
in $\mathbb{P}^3_R$.  In terms of affine coordinates $x=X/Z$ and
$y=Y/Z$, this is the Weierstrass curve with equation
\[
y^2+axy+y=x^3.
\]
Over the extension ring $S=R[D^{-1}]$ with $D=a^3-27$, this
is a (smooth) elliptic curve.

It is more convenient for us to consider an affine neighborhood about
the point $O=[0:1:0]\in \mathbb{P}^3$.  Thus, if we take $u=X/Y$ and
$v=Z/Y$, our curve has the equation
\[
v^2+auv+v=u^3,
\]
where $(u,v)=(0,0)$ is the identity of the elliptic curve.

Let $Q$ be a point on this curve with coordinates $(u(Q),
v(Q))=(d,e)$, and let $\phi\colon C\ra C$ be the map defined in terms
of the group law on the elliptic curve $C$, by $\phi(P)=P-Q$.
Standard considerations allow us to compute that
\begin{align*}
  u(\phi(P)) &= m(-P)^2+am(-P)-d+\tfrac{v}{u^2},
\\
v(\phi(P)) &= m(-P)(u(\phi(P))-d) +e,
\end{align*}
where
\[
m(-P)= \frac{-\frac{v^2}{u^3}-e}{-\frac{v}{u^2}-d}
\]
for $P=(u,v)$. 
Inversion is given by
\begin{align*}
u(-P) &= \frac{v^2}{u^2}+ a\frac{v}{u} -u = -\frac{v}{u^2},
\\
v(-P) &= \frac{v}{u} u(-P) = -\frac{v^2}{u^3}.
\end{align*}

If $Q$ is a point of order $2$, applying the inversion formula shows
that $e=-d^3$ and $d^3-ad-2=0$.  The universal
example $E$ of a subgroup of order $2$ in $C$ is defined over 
\[
S_2 = S[d]/(d^3-ad-2),
\]
and is generated by the point $Q$ with $(u(Q),v(Q))=(d,-d^3)=(d,-ad-2)$.  

Let $C_2$ be the curve with projective equation 
\[
Y^2Z+(a^2+3d-ad^2)XYZ+YZ^2-X^3=0.
\]
In the neighborhood of the basepoint, with $u'=X/Y$ and $v'=Z/Y$, this
has equation
\[
(v')^2+(a^2+3d-ad^2)u'v'+v'=(u')^3.
\]
Let $C_1$ be the curve obtained from $C$ by base change along the
evident inclusion $S\ra S_2$.  There is an isogeny 
\[
\psi\colon C_1\ra C_2
\]
of curves over $S_2$, whose kernel is precisely the rank $2$ subgroup
$E$.
The isogeny can be described easily in terms of $(u,v)$-coordinates:
if $(u(P),v(P))=(u,v)$ in $C_1$ and $(u(\psi(P)),v(\psi(P)))=(u',v')$
in $C_2$, then we have
\[
u' = - u(P)\, u(\phi(P)), \qquad v' = v(P)\, v(\phi(P)),
\]
using the function $\phi$ defined above, where $Q$ is a point of order
$2$.  (This is a kind of ``Lubin
isogeny'' construction.)  

The coordinate $u$ is a uniformizer at the basepoint.  Taking the
formal expansion of $u'$ in terms of $u$, we have
\begin{align*}
u' &= -d\, u+ (ad+3)\, u^2 +(-a^2d-3d^2-2a)\,u^3
+(a^3d+5ad^2+2a^2+6d)\,u^4
\\ & \qquad
+(-a^4d-7a^2d^2-2a^3-16ad-12)\,u^5
\\ & \qquad
+(a^5d+9a^3d^2+2a^4+30a^2d+12d^2+32a)\,u^6 +\cdots,
\end{align*}
and the formal expansion of $v'$ in terms of $u$ is
\begin{align*}
v' &= +(-ad-2)u^3+(2a^2d+3d^2+4a)u^4+(-3a^3d-9ad^2-6a^2-9d)u^5
\\ & \qquad
+(4a^4d+18a^2d^2+8a^3+35ad+23)u^6
\\ & \qquad
+(-5a^5d-30a^3d^2-10a^4-86a^2d-27d^2-84a)u^7
\\ & \qquad
  +(6a^6d+45a^4d^2+12a^5+170a^3d+126ad^2+199a^2+63d)u^8 +\cdots.
\end{align*}
(The easiest way to derive the equation for the curve $C_2$ seems to
be to
calculate the series expansion for $u'$ and $v'$ as above, and then
solve for the Weierstrass equation they satisfy.)

Let 
\[
S_{2,2} = S[d,d']/(d^3-ad-2, (d')^3-(a^2+3d-ad^2)d'-2).
\]
This is a pushout in the category of commutative rings of the diagram
$S_2 \xla{s^*} S \xra{t^*} S_2$, where $s^*\colon S\ra S_2$ is the
usual inclusion, and $t^*\colon S\ra S_2$ sends $a\mapsto
a^2+3d-ad^2$.  
The two ring homomorphims $e_1^*,e_2^*\colon S_2\ra S_{2,2}$ are given
by $e_1^*(a)=a, e_1^*(d)=d$ and $e_2^*(a)=a^2+3d-ad^2, e_2^*(d)=d'$;
they classify the subgroup $E_1$ of 
$C_1$, and $E_2$ of $C_2$, respectively.
The ring $S_{2,2}$ carries the universal example of a
nested chain $(E_1<E_2)$ of subgroups of $C_1$, with $\rank E_1=\rank
E_2/E_1=2$.

Let $S_4$ be defined as the pullback in the following square of
commutative rings.
\[
\xymatrix{
{S_4} \ar[r]  \ar[d]
&{S_{2,2}} \ar[d]_{f^*}
\\
{S} \ar[r]_{s^*}
& {S_2}
}\]
The map $f^*$ sends $d\mapsto d$ and $d'\mapsto a-d^2$; it classifies
the chain of subgroups $(E<C[2])$ in $C$, where $E$ is the universal
example of a rank $2$ subgroup, and $C[2]$ is the subgroup of
$2$-torsion. 
  The universal
example of a subgroup of rank $4$ of $C$ is defined over the ring
$S_4$; the map $S_4\ra S$ classifies $C[2]$.

\section{A Morava $E$-theory spectrum}

Let $\hS=\Z_2\powser{a}$; this is an extension of the ring $S$ of the
previous section.  Let
$\widehat{C}$ denote the formal completion at the identity of the
elliptic curve of the previous section.  This defines a formal group
over $E_0$, which turns out to be a universal deformation for its
reduction to $\F_2=\hS/(2,a)$, a height $2$ formal group.  Let $E$ denote
the even periodic cohomology theory associated to this formal group
over $E_0$.  This is a Morava $E$-theory; it is a strictly commutative
ring spectrum.  

Write $\hS_2 = \hS\otimes_{S}{}_{s^*} S_2$, and similarly for
$\hS_{2,2}$, $\hS_4$, etc.  It is a theorem of Strickland
\cite{strickland-morava-e-theory-of-symmetric} that 
\[
\hS_{p^k} \approx E^0B\Sigma_{p^k}/(\text{transfer}),\qquad k\geq0.
\]

If $F$ is a $K(2)$-local commutative $E$-algebra,
then $\pi_0F$ carries the structure of an amplified $\Gamma$-ring, as
defined above.   The ``total square''
\[
P\colon \pi_0 F \ra \pi_0 F\otimes_{E_0} E^0
B\Sigma_2/(\text{transfer}) \approx \pi_0F \otimes_{E_0}{}_{s^*}\hS_2
\]
can be recovered by the identification
\[ 
\hS_2\approx E^0B\Sigma_2/(\text{transfer}) \approx
\Z_2\powser{a}[d]/(d^3-ad-2), 
\] 
and the formula
\[
P(x) = Q_0(x)+Q_1(x)d+Q_2(x)d^2.
\]
The cartan formula is read off from the formula $P(xy)=P(x)P(y)$.  The
commutation relations are obtained from $P(ax)=P(a)P(x)$, using
\[
P(a)= a^2+3d-ad^2.
\]

Since $\hS_2$ is a free $\hS$-module, there are trace and norm
functions $\hS_2\ra \hS$.  The operators $T$ and $N$ are defined by
$T=\operatorname{trace} P$ and $N=\operatorname{norm} P$.

Let $\mathrm{can}^*\colon \hS_2\ra \hS/(2)\approx \F_2\powser{a}$ be
the map defined by 
$d\mapsto 0$; this classifies the ``canonical subgroup''.  The
operation $P$ has the property that the composite
\[
\pi_0 F \ra \pi_0F \otimes_{\hS}{}_s^* \hS_2 \ra \pi_0F \otimes_{\hS}
{}_{s^*} \hS/(2) \approx \pi_0F/(2)
\]
is precisely the map sending $x\mapsto x^2$; this is the ``frobenius
congruence''.

Applying $P$ twice gives a map
\[
PP\colon \pi_0F \ra 
\pi_0
F\otimes_{E_0}E^0B\Sigma_2/(\text{transfer})_{P}\otimes_{E_0}
E^0B\Sigma_2/(\text{transfer}),
\]
that is, a map
\[
\pi_0F\ra \pi_0F\otimes_{\hS}{}_{s^*} \hS_2{}_{t^*}\otimes_{\hS}{}_{s^*}\hS_2
\approx \pi_0F\otimes_{\hS}{}_{s^*} \hS_{2,2}.
\]
The operation $PP$ factors through $\pi_0F\otimes_{\hS}{}_{s^*}
\hS_4$, where $\hS_4=\lim (\hS_{2,2} \xra{f^*} \hS_2 \xla{s^*} \hS)$; this
is the source of the adem relations. Explicitly, we have 
\begin{align*}
  PP(x) &= P(Q_0x+d\,Q_1x+d^2\,Q_2x)
\\
&= PQ_0x+Pd\,PQ_1x+P(d^2)\,PQ_2x
\\
&= PQ_0x +d'\,PQ_1x + (d')^2\, PQ_2x
\\
&=
Q_0Q_0x+d\,Q_1Q_0x+d^2\,Q_2Q_0x+d'\,Q_0Q_1x+dd'\,Q_1Q_1x+d^2d'\,Q_2Q_1x
\\ &\qquad
+ (d')^2\,Q_0Q_2x + d(d')^2\,Q_1Q_2x + d^2(d')^2\,Q_2Q_2x
\end{align*}
in $\pi_0F\otimes_{E_0}{}_{s^*} \hS_{2,2}$.  Observe that $P(d)=d'$
here, since $P\colon \hS_2\ra \hS_{2,2}$ corresponds to $e_2^*$.  The
projection of the above formula under $f^*\colon \hS_{2,2}\ra \hS_2$ is
\begin{align*}
& Q_0Q_0x +a\,Q_0Q_1x -2\,Q_1Q_1x  +a^2\,Q_0Q_2 x -2a\,Q_1Q_2x +4\,Q_2Q_2
x
\\ &+ (Q_1Q_0x -2\,Q_2Q_1x+2\,Q_0Q_2x)d 
+(Q_2Q_0x -Q_0Q_1x -a\,Q_0Q_2x +2\,Q_1Q_2x)d^2.
\end{align*}
The adem relations are obtained by setting the coefficients of $d$ and
$d^2$ to $0$.  The operation $\Psi$ is what's left over; i.e., it
arises from the
homomorphism $\hS_4\ra \hS$. 

Let $\Gamma'$ denote the algebra of power operations for $E$ as
described in \cite{rezk-morava-e-theory-congruences}.  It has a direct
sum decomposition $\Gamma'\approx \bigoplus \Gamma'[k]$, where the
pieces come from the $E$-homology of $B\Sigma_{p^k}$.  Let
$\Gamma_{\hS} = \hS\otimes_R \Gamma$.  There is a degree preserving
ring homomorphism $\alpha\colon \Gamma_{\hS}\ra \Gamma'$, which is an
isomorphism in degrees $0$ and $1$.  It is easy to see that $\Gamma'$
is generated by $\Gamma'[1]$ using a simple transfer
argument, and thus $\alpha$ is surjective.  We conclude that $\alpha$
is an isomorphism for degree reasons, using the rank calculations of
\cite{strickland-turner-rational-morava-e-theory}.

We have $E^0\mathbb{CP}^\infty\approx \hS\powser{u}$; information about the
action of $\Gamma$ on this ring can be read off from the formal
expansion of $u'$ in terms of $u$ given in \S3.  Thus
\begin{align*}
Q_0(u) &= -3\, u^2 -2a\,u^3
+2a^2\,u^4 +(-2a^3-12)\,u^5 
+(2a^4+32a)\,u^6 +\cdots,
\\
Q_1(u) &= - u+ a\, u^2 -a^2\,u^3
+(a^3+6)\,u^4
+(-a^4-16a)\,u^5
+(a^5+30a^2)\,u^6 +\cdots,
\\
Q_2(u) &=
-3\,u^3
+5a\,u^4
+-7a^2\,u^5
+(9a^3+12)\,u^6 +\cdots.
\end{align*}
The action of $\Gamma$ on $E^0S^2$ is easily read off from this; the
$\Gamma$-module $\omega$ is the kernel of $E^0S^2\ra E^0$.

It is not surprising that the apparatus of ($2$-primary) power
operations can actually be defined over the ring $S=\Z[a,D^{-1}]$,
since all the structure flows from the existence of the elliptic curve $C$
over this ring.  It \emph{is} a bit surprising that all this apparatus
appears to 
be defined over the ring $R=\Z[a]$, where the curve can fail to be
smooth.  However, it appears that all the formulas lift to this
setting, and this is the ground ring I chose to use for the
presentation of \S2.

\newcommand{\noopsort}[1]{} \newcommand{\printfirst}[2]{#1}
  \newcommand{\singleletter}[1]{#1} \newcommand{\switchargs}[2]{#2#1}
\providecommand{\bysame}{\leavevmode\hbox to3em{\hrulefill}\thinspace}
\providecommand{\MR}{\relax\ifhmode\unskip\space\fi MR }
\providecommand{\MRhref}[2]{%
  \href{http://www.ams.org/mathscinet-getitem?mr=#1}{#2}
}
\providecommand{\href}[2]{#2}

\end{document}